\documentclass{amsart}

\usepackage[english] {babel}
\usepackage{amsmath,amssymb,amscd}
\usepackage{amsthm}
\usepackage{verbatim}
\usepackage{ifthen}
\usepackage[all]{xy}
\usepackage{bookmark,hyperref}
\hypersetup{
     colorlinks=true,
     linkcolor=blue,
     filecolor=blue,
     citecolor = blue,
     urlcolor=cyan,
     }

\theoremstyle{definition}
\newtheorem*{defi}{Definition}

\theoremstyle{plain}
\newtheorem*{thmsn}{Theorem}
\newtheorem{teo}{Theorem}

\theoremstyle{remark}

\makeatletter
\@namedef{subjclassname@2020}{%
  \textup{2020} Mathematics Subject Classification}
\makeatother

\emergencystretch=\maxdimen
\begin{document}

\title[Partially hyperbolic motions in the N-body problem]{Existence of partially hyperbolic motions in the N-body problem}

\author{J. M. Burgos}

\address{Departamento de Matem\'aticas, Centro de Investigaci\'on y de Estudios Avanzados, Av. Instituto Polit\'ecnico Nacional 2508, Col. San Pedro Zacatenco, C.P. 07360 Ciudad de M\'exico, M\'exico.}

\email{burgos@math.cinvestav.mx}

\begin{abstract}
In the context of the Newtonian $N$-body problem, we prove the existence of a partially hyperbolic motion with prescribed positive energy and any initial collisionless configuration. Moreover, it is a free time minimizer of the respective supercritical Newtonian action or equivalently a geodesic ray for the respective Jacobi-Maupertuis metric.
\end{abstract}



\maketitle



\section{Introduction}

In this short note we deduce, from the recent theorem by Maderna and Venturelli \cite{MV} showing the existence of hyperbolic motions with arbitrary limit shape in the classical $N$-body problem, the following corollary:

\begin{teo}
Within any positive energy level and starting at any given collisionless configuration, provided that the underlying space has dimension at least two, there is a partially hyperbolic motion.
\end{teo}

This is a motion defined on some closed half real line in the future that dynamically evolve as clusters separating linearly in time such that the mutual distance between two bodies in the same cluster grows like $t^{2/3}$ and there is at least two clusters with one of them having more than one body. We prove the theorem by showing the existence of a free time minimizer partially hyperbolic motion of the $h$ supercritical Newtonian action.

Consider an Euclidean space $E$ with $\dim E\geq 1$ and define the configuration space as the cartesian product $E^{N}$ with the \textit{mass inner product}
$$\langle x,y\rangle= \sum_{i=1}^{N}\,m_i\langle x_i, y_i\rangle$$
where $m_1,\ldots m_N>0$ denote the masses of the respective bodies. With respect to this geometry, Newton's motion equation acquires the simpler form
$$\ddot{x}= \nabla U(x),\qquad x\in \Omega$$
where $U$ is the Newtonian potential and $\Omega$ denotes the open dense set of collisionless configurations
$$\Omega= \{(x_1,\ldots x_N) \mid x_i\neq x_j\ \ {\rm if}\ \ i\neq j\}.$$

Given a pair of configurations $x,\ y$ and $\tau>0$, we denote the space of absolute continuous joining $x$ with $y$ in time $\tau$ by $\mathcal{C}(x,y,\tau)$. Define the fix and free time \textit{critical potentials} respectively by
$$\phi(x,y,\tau)= \inf\{A_L(\gamma) \mid \gamma\in \mathcal{C}(x,y,\tau)\},$$
$$\phi(x,y)= \inf\{\phi(x,y,\tau) \mid \tau>0\}$$
where $A_L$ denotes the action of the Newtonian $N$-body Lagrangian $L$. Every minimizer of the action among this space is a critical curve hence it is a solution of the Newton's equation whenever it does not have a collision. However, as it was earlier noticed by Poincar\'e in \cite{Poincare}, there are finite action curves with isolated collisions hence a minimizer orbit could a priori not be a true motion.

A breakthrough in the theory was given by the following Theorem which unlocks the use of variational methods in the Newtonian $N$-body problem. The main idea is due to Marchal in \cite{Marchal} and complete proofs were given by Chenciner in \cite{Chenciner} and Ferrario and Terracini in (\cite{Terracini}, Corollary 10.6).

\begin{thmsn}
Provided that $\dim E\geq 2$, minimizers of the fixed ends problem are free of interior collisions. That is to say, if $\gamma:[a,b]\to E^{N}$ is an absolute continuous curve such that $A_L(\gamma)= \phi(\gamma(a), \gamma(b), b-a)$, then $\gamma(t)\in \Omega$ for every $t$ in $(a,b)$.
\end{thmsn}

Define the \textit{supercritical action potential} at the energy level $h$ by
$$\phi_h(x,y)=\inf \{A_L(\gamma)+h\tau \ |\ \gamma\in \mathcal{C}(x,y,\tau),\ \tau>0\}.$$
Because the Ma\~n\'e's critical value of the Newtonian $N$-body Lagrangian is zero, the previous definition only makes sense for $h\geq 0$.

\begin{defi}
An absolute continuous curve $\gamma$ defined on some interval $I$ is an \mbox{$h$ \textit{free time minimizer}} if
$$A_{L+h}\left(\gamma|_{[t_1,t_2]}\right) = \phi_h\left(\gamma(t_1), \gamma(t_2)\right)$$
for every $t_1$ and $t_2$ in the interval $I$ such that $t_1<t_2$.
\end{defi}

In particular by Marchal's Theorem, if $\dim E\geq 2$ and $h\geq 0$, then $h$ free time minimizer curves defined on the interval $[a,b]$ are true motions on $(a,b)$ with energy $h$.

Equivalently, the main result can be expressed geometrically in terms of the Jacobi-Maupertuis Riemannian metric
$$j_h=2(h+U)\,g|_\Omega$$
at the energy level $h\geq 0$ over $\Omega$ where $g$ denotes the flat metric induced on $E^N$ by the mass inner product. Trajectories of a fixed energy level $h$ are geodesics of this metric.

\begin{defi}
\label{def:GR}
Consider $h\geq 0$. A curve $x:[t_0,+\infty)\to\Omega$
is a \emph{geodesic ray} for the metric $j_h$ whenever
\emph{its arclength parametrization} is an isometric embedding of
the half line $[t_0,+\infty)$ within the Riemannian space $(\Omega,j_h)$.
\end{defi}

It is well known the equivalence between the variational free time minimizer property of $A_{L+h}$ and the geometrical geodesic ray property for the Jacobi-Mapertuis metric $j_h$. In particular, we prove

\begin{teo}
Within any energy level $h>0$ and starting at any given collisionless configuration, provided that $\dim E\geq 2$, there is a partially hyperbolic geodesic ray for the Jacobi-Mapertuis metric.
\end{teo}

Denote by $\mathcal{C}_{t_0}$ the space of non superhyperbolic motions defined on $[t_0,+\infty)$ with the topology of $C^{1}$ convergence over compact sets. Consider the map associating to every non superhyperbolic motion its final configuration
$$C:\mathcal{C}_{t_0}\rightarrow E^{N},\qquad \left( a\,t+o(t)\right)\mapsto a.$$
With the usual topology on the space of configurations, since the work of Chazy \cite{Ch2} and called by him as \textit{continuit\'e de l'instabilit\'e}, it is known that this map is continuous on the subspace of hyperbolic motions, i.e. those whose final configuration is collisionless (See also Lemma 4.1 in \cite{MV} for a modern proof).

We conjecture that in general, with an arbitrary number of bodies with arbitrary masses, this map cannot be continuously extended to the whole space $\mathcal{C}_{t_0}$ and this is the main difficulty in the control of the final configuration of the partially hyperbolic motion in our result. Moreover, we conjecture that even with a converging sequence in $\mathcal{C}_{t_0}$ to some partially hyperbolic motion, the omega limit of the image of this sequence by the map $C$ has multiple points. Interpreting the existence of clusters with multiple bodies as resonances, this phenomenon would be the analog of \textit{Takens chaos} in the $N$-body problem.

\section{Proof}

\begin{proof}
Consider the \textit{center of mass} linear map $G:E^{N}\rightarrow E$ that is the composition of the orthogonal projection with respect to the mass inner product onto the diagonal and the canonical isomorphism of the diagonal with $E$. Concretely,
$$G(x)= \sum_{i=1}^{N}\,m_i\,x_i,\qquad x=(x_1,\ldots, x_N)\in E^{N}.$$

Let $h>0$ and consider a sequence $(a_n)$ of configurations in $\Omega$ such that
$$||a_n||^{2}/2= h, \qquad G(a_n)=0$$
for every $n$ and converging to some $b$ in $E^{N}-\Omega$. Note that necessarily by continuity
$$||b||^{2}/2=h>0,\qquad G(b)=0$$
and in particular $b$ is not zero.

Let $x_0\in \Omega$. By (Theorems 3.2 and 3.4 in \cite{MV}), for every $n$ there is an $h$ free time minimizer hyperbolic motion $\gamma_n:[0,+\infty)\rightarrow \Omega$ such that
$$\gamma_n(t)=a_n t+o_n(t),\quad\quad \gamma_n(0)=x_0.$$
The initial velocities $v_n=\dot{\gamma}(0)$ lie in a sphere and taking a subsequence if necessary, we may assume that there is $v$ in the sphere such that $(v_n)$ converges to it.

Let $\zeta:[0, \omega_+)\rightarrow\Omega$ be a solution with maximal $\omega_+$ such that $\zeta(0)=x_0$ and $\dot{\zeta}(0)=v$. It is clear that $\zeta$ has energy $h$ and $\omega_+>0$ for $x_0$ is in $\Omega$. We will prove that $\omega_+=+\infty$ and  $\zeta$ is a partially hyperbolic free time minimizer motion.

Consider a sequence $(\lambda_n)$ such that $o_n(\lambda_n)= o(\lambda_n)$ and define $p_n= \gamma_n(\lambda_n)$. Therefore, the sequence $(p_n)$ verifies the following asymptotics:
$$p_n= a_n\lambda_n+o_n(\lambda_n)= b\lambda_n+(a_n-b)\lambda_n+o_n(\lambda_n)= b\lambda_n+o(\lambda_n).$$

For every configuration $p$, define the continuous function $u_{p}$ by the formula $u_p(x)= \phi_h(x,p)$ and note that it is a Hamilton-Jacobi viscosity subsolution for it is dominated by $L+h$:
$$u_p(x)-u_p(y)= \phi_h(x,p)-\phi_h(y,p)\leq \phi_h(x,y).$$
Moreover, it is clear that every $h$ free time minimizer curve $\xi:[a,b]\rightarrow E^{N}$ such that $p=\xi(b)$ is an $h$ calibrating curve of $u_p$ for it is trivially verified that
$$u_p(\xi(a))-u_p(\xi(b))= A_{L+h}(\xi).$$
The only point at which $u_p$ fails to be a global viscosity solution is $p$.

The set of Hamilton-Jacobi viscosity subsolutions vanishing at the origin is compact (Corollary 2.12 in \cite{MV}) with respect to the compact-open topology hence, taking a subsequence if necessary, we may assume that the sequence $(p_n)$ defines the Hamilton-Jacobi viscosity subsolution
\begin{equation}\label{eq0}
u(x)= \lim_{n\to +\infty} \left(u_{p_n}(x)-u_{p_n}(0)\right).
\end{equation}
This is a \textit{horofunction directed by $b$} and it is a global Hamilton-Jacobi viscosity solution (Theorem 3.1 in \cite{MV}).

Suppose that $\omega_+$ is finite. There is a natural $n_0$ such that $\lambda_n\geq \omega_+$ for all $n\geq n_0$. By the continuity with respect to the parameters, considering the subsequence $(\gamma_n)_{n\geq n_0}$ we have that $\gamma_n\to \zeta$ and $\dot{\gamma}_n\to \dot{\zeta}$ uniformly over compact sets of $[0, \omega_+)$. Then,
\begin{equation}\label{eq1}
\lim_{n\to +\infty} A_{L+h}(\gamma_n|_{[0,\tau]}) = A_{L+h}(\zeta|_{[0,\tau]})
\end{equation}
for every $0< \tau<\omega_+$. Because each $\gamma_n$ is an $h$ free time minimizer motion,
\begin{equation}\label{eq2}
A_{L+h}(\gamma_n|_{[0,\tau]})= u_{p_n}(x_0) - u_{p_n}(\gamma_n(\tau))
\end{equation}
for every $0< \tau<\omega_+$ and every $n\geq n_0$. Because the convergence in \eqref{eq0} is uniform over compact sets, from equations \eqref{eq1} and \eqref{eq2} we have
$$u(x_0) - u(\zeta(\tau))= A_{L+h}(\zeta|_{[0,\tau]})$$
for every $0< \tau<\omega_+$ hence $\zeta$ is an $h$ calibrating curve of $u$.

Let $0<t_*<\omega_+$. There is an $h$ calibrating curve (Theorem 3.2 in \cite{MV})
$$\gamma:[0,+\infty)\rightarrow E^{N}$$
of the horofunction $u$ such that $\gamma(0)=\zeta(t_*)$. Because the concatenation of calibrating curves is also calibrating, the concatenation of $\zeta|_{[0,t_*]}$ with $\gamma$ is also an $h$ calibrating curve hence an $h$ free time minimizer curve. By Marchal's Theorem, this concatenation is a true motion verifying $\dot{\gamma}(0)=\dot{\zeta}(t_*)$ and it is a proper extension of $\zeta$. This is clearly a contradiction for $\omega_+$ was maximal and we conclude that $\omega_+= +\infty$.

Because $\zeta:[0,+\infty)\rightarrow E^{N}$ is an $h$ calibrating curve, it is an $h$ free time minimizer motion and in particular it is not superhyperbolic (\cite{MV}, argument in page 31). By (Theorem 1 in \cite{MS}), there is a configuration $b'$ in $E^{N}$ such that
$$\zeta(t)= b't+o(t).$$

If $b'$ is in $\Omega$, then by Chazy's \textit{continuit\'e de l'instabilit\'e} we have $a_n\to b'$ hence, by the limit uniqueness, $b=b'$ so $b$ is in $\Omega$ as well, a contradiction.

Because $G(a_n)=0$ for every $n$, every $\gamma_n$ has fix center of mass $G(x_0)$ and by continuity the same occurs with $\zeta$. In particular, $G(b')=0$. If $b'$ is a total collision, then it must be zero for its center of mass is zero but this is a contradiction for $\zeta$ has energy $h>0$.

We conclude that $b'$ neither is a total collision nor it is in $\Omega$. We have proved that $\zeta$ is an $h$ free time minimizer partially hyperbolic motion.
\end{proof}

\section*{Acknowledgments}

The author is grateful with \textit{Consejo Nacional de Ciencia y Tecnolog\'ia} for its \textit{C\'atedras Conacyt} program. He is also grateful with Professor Ezequiel Maderna for helpful conversations.

\end{document}